# Nelson Goodman's *qualia* and the Discrete Cat Mapping

Godofredo Iommi Amunátegui
Instituto de Física
Pontificia Universidad Católica de Valparaíso
giommi@ucv.cl

Abstract

In *The Structure of Appearance* Nelson Goodman aims at describing the world as composed of individuals. He builds his system on phenomenal qualities or *qualia*. We propose a mathematical model for the *qualia* based on a dynamical system, the so-called "Discrete Cat Mapping".

Key Words: Nominalism, Qualia, Dynamical Systems.

In *The Structure of Appearance*[1], Nelson Goodman pose, expose and compose the core of his thought[2]. Hereafter we intend, for the most part, to establish a link between the basic units or *qualia* of Goodman's system and the chaotic mapping called Arnold Cat Map[3]. Sections I and II are dedicated to put forward an overview of Goodman's opus that points out the traits relevant for our purpose. Section III treats the Discrete Cat Map and Section IV deals, properly, with the proposed interpretation. In these pages we take, in general, the liberty of repeating the chosen passages of Goodman's works without quotation marks. As a bonus, this procedure allows us to avoid awkward paraphrases.

I. Nominalism

Nominalism describes the world as composed of individuals. To explain nominalism we need to explain not what individuals are but rather what constitutes describing the world as composed of them. So to describe the world is to describe it as made up of entities no two of which break down into exactly the same entities. Nominalism does not involve excluding abstract entities but requires only that whatever is admitted as an entity at all be construe as an individual[4].

Goodman's nominalism "is a positions that [ … ] renounces classes and other non-individuals. The position that *also* accepts non-individuals he calls "Platonism" "[5].

Let us suppose, for example, that a nominalist and a Platonist start with the same minimal, atomic elements for their systems. The nominalist admits also all wholes or individuals sums comprised of these, and so has a universe of a finite number of entities. The Platonist admits no sums of atoms but admits all classes of them, and he further admits all classes of classes of atoms, and so on *ad infinitum.* He gets all these extra entities out of his original atoms. Goodman's dictum is "No distinction of entities without distinction of content".

## II. *Qualia*
### A. Definitions and Determinations

To say that a thing *looks* green is to make a statement concerning a presented quality, a color quality of some presentation of the thing, while to say that a thing *is* green is to make a more complex statement concerning the color qualities exhibited by various presentations of the thing. The color names are thus used in two different ways in ordinary language: in the one case for presented characters, hereafter called *qualia*; in the other, for properties of things[6].

Though presentations are momentary and unrecallable, they are nevertheless comparable in that they contain repeatable and recognizable *qualia*. Any judgment that a *quale* of one presentation is the same as a quale of another is open to pertinent criticism that may cause it to be abandoned. If we divide the stream of experience into its smallest concrete parts and then go on to divide these concreta into sense *qualia*, we arrive at entities suitable as atoms for a realistic system. A visual concretum might be divided, for example, into three constituent parts: a time, a visual-field place, and a color[7].

A visual concretum is already a spatially smallest discernible particle of phenomena and the further analysis into its three component *qualia* leaves its space undivided. This analysis consists simply of distinguishing the place from both the time and color that together with the place, make up the concretum.

Acceptance of *qualia* as individuals is not inconsistent with the refusal to countenance classes and other nonindividuals. Nominalism excludes all except individuals but does not decide what individuals are. To treat *qualia* as atoms for a system as not to suppose that they are the units in which experience is originally given. Even should it be shown that certain units of experience are in some sense apprehended prior to any operation of analysis or synthesis, this would not at all preclude the selection of other units as atoms, for a constructional system is not necessarily intended as an epistemological history.

Location in the visual field may change when physical position is constant, or remain constant when physical position changes. A grammatical convention requires us to say of a color *h* and a place *k* that *h* occurs at *k* and not that *k* occurs at *h*; but this can be construed as a condensed way of saying that the two occur together[8]. By saying that two *qualia* are together, we intend to make an assertion having the same nontemporal character as a mathematical statement: a color is as at a place if the color was, is now or will be at that place and a color occurs at a time even though the time be past or future. Thus such expressions is "occurs at" or "is together with", like "intersects" in mathematics, are to be taken as without effective tense[9].

## B. Topology

The specific *qualia* that an individual contains fix its size and shape with respect to each category of *qualia*. If we can count and are familiar with the order of *qualia* in a given category, then we can determine the size and shape of an individual *x* in that respect if we know exactly what *qualia* of that kind are contained in *x*.

For the arrangement of a given set of *qualia* is not variable. *Qualia* cannot be "moved around"; each has –we might almost say each *is*– a fixed position in the array of the category to which it belongs. Terms for color, place, time, etc., may be called initial quality terms, while terms for shape and size are called derivative quality terms. Initial quality terms are construed as names of individuals while derivative quality terms are construed as syncategorematic. The central problem is to construct, for each category of *qualia*, a map that will assign to each *qualia* in the category a unique position and that will represent relative likeness of *qualia* by relative nearness in position[10].

We never actually compare more than two or three very similar *qualia* at any moment. Any adequate map of a category must thus be constructed on the basis of comparison made at many different times. Identifications of *qualia* from one moment to another are in a sense decrees but are made and corrected with a view to the difficulties than result from an ill-considered set of decrees. Identifications of physical objects and identification of *qualia* must of be co-ordinated with one another. The actual process of mapping is affected by practical limitations. We begin with observations of what we hope may prove to be a typical set of very similar *qualia*, seek to determine the pattern exhibited, and then by interpolation, extrapolation, and conjecture construct a tentative complete map that may later be improved and corrected[11].

A familiar "kind" or "category" of *qualia* consists of all colors, or all places, or all sounds, etc. When we compare two concreta with respect, say, to position or color, we are comparing a *quale x* that is part of one concretum with a *quale y* that is a part of the other and that belongs to the same category as *x*. Given any two *qualia* belonging to the same category, we can trace a path from one to the other by a series of steps, each to a *quale* matching the preceding one. When two *qualia* belong to different categories there is no such path joining them. We cannot go from a color to a time by such a series of steps. This suggests a way of defining categories: two *qualia* are joined by a path-called an M-path-just in case the ancestral of "M" applies between them. An individual is a clan if it cannot be divided into two parts such that no *quale* in one matches any *quale* in the other. A single *quale* constitutes such a clan, as does the sum of any group of *qualia* that match each other, the sum of two non-matching *qualia* and a third that matches both, and so on. A clan contains an M-path between each two of its *quale* parts.

A *quale-category* is a most comprehensive clan. A category contains every *quale* that matches any part of it. A clan that fails of being a category lacks one or more *qualia* of the category containing the clan. Any sense realm such as that of the visual qualities or that of the auditory qualities, is defined as a sum of certain categories. A realm as a general term may be defined by enumeration of the several realms. The sum of all the *qualia* that match a given *quale* is called a manor. Every *quale* is a part of its own manor and of the manor of any matching *quale*. As an illustration let us consider the positive natural numbers and the following "rule": the *quale* designated by any number matches the *qualia* designated by numbers not more than 4 places removed from that number (we define this a linear M-span of n=4):

(i) 1 2 3 4 5 6 7 8 9 10 11 …
(ii) Manor of 1 = 1 + 2 + 3 + 4 + 5
     Manor of 2 = 1 + 2 + 3 + 4 + 5 + 6
     Manor of 3 = 1 + 2 + 3 + 4 + 5 + 6 + 7

Note that the rule is arbitrary to some degree for we do not have before us a fully articulated presystematic order against which the details of our constructed order can be checked.

The rule is no part of system. It is an extrasystematic statement of a relationship between matching and order. The calculus is designed to apply to finite sets of elements. This fact explains why some of the most familiar theorems concerning order in the continuum, such that there is an element between each two distinct elements, do not hold here. At the same time, the calculus must be adequate for dealing with arrays of any degree of complication. We can by no means take it for granted that all categories of *qualia* are linear arrays[12].

C. <u>Time</u>

Phenomenal time is merely one among the general categories of *qualia*. Every *quale* that is not itself a time occurs at some time; hence every concretum, no matter of what sense realm, contains a time. Consider a colored patch that appears in the visual field and stays for a while. The total presentation comprises the colors, places and times involved. It is temporally as well as spatially divisible; and its identity over different times, like its identity over different places, is the identity of a totality of diverse parts. In contrast, the sum of color-spots that constitutes the colored patch itself contains no times. Neither it nor any of its parts is temporally divisible. It retains its strict numerical identity throughout the period in question. Such occurrence through a period is quite a different matter from the duration of a thing or event; and we had better to observe the distinction by saying that the patch persists through the period. Under these definitions, an individual that endures or persists through a period endures or persists through every part of that period. The endurance or persistence is continuous or discontinuous according as the period is. Thus although a color *quale* occurs at times and persists through periods, it is nevertheless literally out of time. It is eternal. That is not to say that it is everpersisting, for nothing occurs at all times. Nor is an eternal individual everlasting; for an everlasting or ever-enduring individual is one that contains all times. Observe that the eternity of an individual is no bar to its occurrence at some times or its failure to occur at others. Only what is eternal is with a time[13].

## III. The Discret Arnold Cat Mapping

Arnold and Avez develop a study of the classical systems with strongly stochastic properties, the so-called C-systems. As an example they consider a continuous two-dimensional invertible chaos map[14]. When the variables *x*,*y* belong to [0,1,…, N-1] the map is discretized. Through this exposition the mathematical details are reduced to the bare minimum consistent with clarity[15].

Let (*x*,*y*) denote a point in the unit square. The mapping takes (*x*,*y*) to the new point

$$\begin{pmatrix} x' \\ y' \end{pmatrix} = \begin{pmatrix} 1 & 1 \\ 1 & 2 \end{pmatrix} \begin{pmatrix} x \\ y \end{pmatrix} (\mathrm{mod}\, 1).$$

The mapping preserves area[16]. The computer is a convenient device for demonstrating mappings where the screen serves as a two-dimensional lattice of points (pixels). Consider a square lattice of points and denote the points by (*x*,*y*). Their discretized values with the operations of addition and multiplication performed (mod N), may be expressed by the mapping

$$\begin{pmatrix} x' \\ y' \end{pmatrix} = \begin{pmatrix} 1 & 1 \\ 1 & 2 \end{pmatrix} \begin{pmatrix} x \\ y \end{pmatrix} (\mathrm{mod}\, N).$$

N is selected so as to make ample use of the capability of the screen. Starting with the initial configuration (a cat's image, for instance), "snapshots" of the images at different times under the mapping are currently given in the literature[17]. By placing a picture of a cat in the unit square and then displaying several subsequent images resulting from the discrete flow, a mixing property appear: the images show that the cat tends to become "smeared" over the unit square[18]. The initial configuration must eventually return, since there are $2^{N \times N}$ possible configurations of the N×N pixels (each pixel is "on" or "off"). How long does it take an N dimension picture to return to its original image? Dyson and Falk show[19] that the periodicity $m_N$ of the mapping has an upper bound m* which may be evaluated in a rather direct manner:
   (i) Consider a positive integer *N*>1 and write *N* in terms of its prime factors *p* and *q*:

$$N = \left(\prod_{p/N} p^{\alpha}\right)\left(\prod_{q/N} q^{\beta}\right) 5^{\gamma} 2^{\delta}$$

Where *p/N* means "p divides N".

(ii) The upper bound m* results from the formula:

$$2m^* = LCM\left[(p-1)p^{\alpha-1}, 2(q+1)q^{\beta-1}, 2(10)5^{\gamma-1}, (3)2^{\varepsilon}\right]$$

where $\varepsilon = M(\delta-1,1)$; LCM means "compute the least common multiple". Remark that each factor of N written in terms of its factors has a corresponding term in the LCM.

Example: $N = 300 = 2^2 \cdot 3 \cdot 5^2$

$5^2 \to 2(10)5 = 100$

$3(2^{\varepsilon}) \to 3 \cdot 2 = 6$ ($\varepsilon$= Max(2-1,1)=1)

Thus 2m* = LCM (100,6) = 600 and m* = 300.

The upper bound of the period is 300. It must be noted that:

(a) The period of the Discrete Cat Map does not always become greater with an increasing modulo (i.e., for N = 300 and N = 150, m* has the same value)[20].

(b) The number $2^{N \times X}$ is large and it is surprising to see that the configuration returns after only 24 iterations for N = 161 or after 15 iterations for N = 124. In general, the relationship between the size of the image and how many iterations it takes to return to its original form appears to be random.

# IV. Towards a possible interpretation

In this section we intend to establish a relation between some aspects of Goodman's system and the Discrete Cat Map. As a proviso, we propose the following table of equivalences:

| *quale, qualia* | pixel, pixels |
|---|---|
| Two-dimensional concretum | Configuration on the screen |
| Map | Arnold discrete map |

The problem is to construct, for each category of *qualia*, a map that will assign to each *quale* in the category a unique position and that will represent likeness of *qualia* by relative nearness in position. Topology and time determine the boundaries of the problem and some of their features may be treated in terms of the discrete map model. Let us focus on some of them:

(1) Every *quale-pixel* occupies a point that pertains to the configuration of the N×N *qualia-pixels*.

- Identifications of configurations and identifications of *qualia* must be coordinate with one another. The discrete map defines such a coordination.
- Every adequate map of a concretum must be constructed on the basis of comparison made at many different times. The discrete map determines the evolution of each *quale-pixel*.

(2) The image –of a cat, for a change– tends to become "smeared" on the N×N plane, i.e., on the screen. Hence a *quale's* position varies at random between two period-stances. At any time a *quale* is "somewhere" on the screen. According to Goodman's conception *qualia* cannot be "moved around". Each *quale* has a fixed positition in the array of the category to which it belongs. We seem to indulge in a gross approximation when indicating a *quale's* position by means of the adverb "somewhere". But this term points to a mathematical fact which results from the discrete mapping we are dealing with. Moreover Goodman also indicates that identifications of *qualia* from one moment to another are, in a sense, decrees.

(3) Godman's calculus is designed to apply to finite sets of elements and it must be adequate for dealing with arrays of any degree of complication. In a discrete map, for an image of any degree of complexity the number $N$ of *qualia-pixels* is finite.

(4) "The assertion "two *qualia* are together" has the same nontemporal character as a mathematical statement". In the discrete map model, after m* iterations, the same image reappears on the screen and the *qualia* relative positions –in the initial and final states– are the same, i.e., they are constant. What happens during these states? In (2) we have outlined a plausible answer to such a question.

(5) "Phenomenal time is merely one among the general categories of *qualia*". But on the screen the successive configurations themselves ("snapshots") do not display *time-qualia*, although time pertains, properly, to the discrete map definition. At his point, a subtle conceptual turn is sharply unveiled by Goodman: "Every *quale* that is not itself a time occurs at some time". A remarkable property of the model corresponds, in our opinion, to Goodman's statement: the period of the mapping depends on the number $N$ of *qualia-pixels*. Hence, in a way, every *quale* is related to some time.

(6) "The colored patch contains no times. It retains its strict numerical identity through periods, it is out of time". Iterations, periods and *qualia* numerical identity conform structural traits of the map model. Moreover the map is a mathematical conception and, by its very nature, "is out of time". So, for instance, a *color-quale* is at a place if it was, is now or will be at that place and occurs at a time even though the time be past or future.

(7) The sum of *qualia* that match a given *quale* –such is Goodman's definition of a "manor"– is no part of the system. It is an extrasystematic statement between matching and order. The discrete mapping takes the point $(x,y)$ to $(x',y')$ by restricting the variables to the set $[0,1, … N\text{-}1]$.

Suppose that this set of integers arranged according to their natural increasing values fixes a system, and let us introduce two equivalences: (a) between this set and Goodman's system, and (b) between the extrasystematic rule and the discrete mapping.

We remark that the resulting order of the *qualia-pixels* neither represents nor is an intrinsic structure of a given set of integers. Hence, the rule is arbitrary for there is no presystematic order against which the constructed order can be checked.

Our choice of these seven illustrations of the proposed relationship is, somewhat, arbitrary and by no means exhaustive. Moreover each one of them is open to a critical analysis. Nevertheless at least this list of "simple samples" clarifies the intuition which sustains this work.

## V. Conclusion

In order to close the article two appreciations of Goodman appear to be pertinent:
- "The actual process of mapping is naturally affected by practical limitations […] we seldom if ever have the opportunity of starting from a comprehensive set of observations […]. Rather we begin with observations of what we hope may prove to be a typical set of very similar *qualia*, seek to determine the pattern exhibited, and then by interpolation, extrapolation, and conjecture construct a tentative map that may later be improved and corrected"[21].
- "What I have offered here constitutes only the beginning of a theory of nonlinear order. The best hope for further progress seems to me to lie in finding ways of making greater use of algebraic techniques"[22].

Although the analogy displayed through these pages may suffer a number of short comings, it conforms to a main project: to improve the map employing algebraic techniques.

## Acknowledgments

Our interest on Arnold's work is mainly due to Godofredo Iommi Echeverría's sharp views on dynamical systems. This research has been partially supported by Fondecyt (Proyect Nº 1120019).

# References

[1] Nelson Goodman, *The Structure of Appearance*, The Bobbs-Merrill Compo., second edition, 1966 (S.A.).

[2] "Alas, his most important work is also his most complicated, which is why it is so often ignored". This "diagnosis" is formulated by D. Cohnitz and M. Rossberg in *Nelson Goodman*, McGill-Queen's University Press, 2006, p. 139 (C-R). The authors emphasize that a *Study of Qualities* (Goodman's Dissertation), first published in revised form as *The Structure of Appearance*" is best understood as situated within a certain branch of analytic philosophy, departing from Frege's conception of ideal-language philosophy". They indicate that "Rudolf Carnap is his predecessor in this tradition" (C-R, p. 99). Goodman scrutinizes Carnap's *Der logische Aufbau der Welt: Scheinprobleme in der Philosophie*, Hamburg: Meiner, 1961. (first edition 1928). It seems convenient to recall his outline of Carnap's *Aufbau*: "He first decides upon his basic units or "ground elements". Then he describes the method by which he will seek to construct qualities, and with the requirements of that method in mind then selects his primitive relation. Because he chooses a relation that will also provide him with the means for ordering qualities, he finds it necessary to define quality classes by a somewhat, indirect route. That accomplished, he discusses the question how such qualities may be classified into the various sense realms: auditory, visual, tactual, etc. Concentrating next upon the visual qualities, he deals with the problem of separating and distinguishing the spatial qualities, or locations, from the color qualities; this is not another subclassification of the visual qualities but […] a further abstraction from the more concrete color-spots he first abstracted from the ground elements. He deals then with the matter of ordering color and visual places […]". (S.A., Chapter five). Goodman espouses, in part, an analogous strategy.

[3] V.I. Arnol'd and A. Avez, *Ergodic Problems in Classical Mechanics,* translated from the French by A. Avez, W.A. Benjamin Inc., 1968, pp. 53-59. [A-A].

[4] Nelson Goodman, *Problems and Projects*, The Bobbs-Morrill Company, Inc., 1972, pp. 158-163.

[5] C-R, p. 137.

[6] S.A., p. 130. In a footnote Goodman indicates that the usage of the term *quale is taken from C.I. Lewis "Mind and World order", Scribner's Sons, N.Y. 1929.*

[7] S.A., pp. 132-135; p. 189.

[8] S.A., pp. 191-195.

[9] S.A., p. 203

[10] S.A., pp. 260-261; p. 268.

[11] S.A., pp. 278-280.

[12] S.A., p. 285-297; pp. 300-301; pp. 323-328. Goodman analyzes some types of nonlinear arrays: square-cell network whose maximum manor has $1+2n(n+1)$ elements; triangular- cell networks whose maximum manor has $1+3n(n+1)$ elements; cubical-cell networks whose maximum manor has $1+2n + 2n(n+1)(2n+1)/3$ elements

[13] S.A., pp. 355-358. The paragraph's very last sentence might be read *cum grano salis*: "Theologians have perhaps overlooked something here"..

[14] A-A, pp. 53-54.

[15] We follow, closely, the beautiful article of Freeman J. Dyson and Harold Falk, "Period of a Discrete Cat Mapping", The American Mathematical Monthly, Vol. 99, 7, 1992, pp. 603-614 (D-F).

[16] The sentence –"A mapping having the above mathematical properties and connections with statistical mechanics has an "intellectual domain of attraction", and we were drawn in. This paper documents our pleasant experience"- shows, in my opinion, the Dyson-touch.

[17] D-F, p. 605. In fact, it suffices to open a link …

[18] D-F, p. 604.

[19] For details, see D-F, pp. 607-608.

[20] J. Bao and Q. Yang, "Period of the discrete Arnold cat map and general cat map", Nonlinear Dynamics, 70, 2012, pp. 1365-1375; A table for the period of the discrete cat map from $N=2$ to $N=100$ may be consulted at p. 1369. Let us point out that the irregularity of the behavior of iterates of points has been studied (T-Y Li and James A. Yorke, "Period Three Implies Chaos", The American Mathematical Monthly, Vol. 82, 10, 1975, pp. 985-992).

[21] S.A., p. 280.

[22] S.A., p. 330.